\documentclass{fundam}
\usepackage{url} 
\usepackage{graphicx}

\makeatletter
\def\ps@pprintTitle{%
 \let\@oddhead\@empty
 \let\@evenhead\@empty
 \def\@oddfoot{\centerline{\thepage}}%
 \let\@evenfoot\@oddfoot}
\makeatother

\usepackage{mathrsfs}
\usepackage{algpseudocode}
\usepackage{amssymb, xfrac}

\usepackage{color,soul}
\usepackage[shortlabels]{enumitem}
\usepackage{todonotes}
\usepackage{amsmath}
\usepackage[ruled,vlined]{algorithm2e}
\usepackage{soul}
\newtheorem{thm}{Theorem} 
\newtheorem{lem}[thm]{Lemma}

\newtheorem{obs}[thm]{Observation}
\newtheorem{rem}[thm]{Remark}
\newtheorem{prop}[thm]{Proposition}

\newtheorem{sub-alg}[thm]{Sub-Algorithm}

\newtheorem{deff}[thm]{Definition}

\newenvironment{myproof}[1][Proof:]%
{\trivlist\PRstyle\item[]{\bfseries#1}\newline}%
{\QED\endtrivlist}

\newcommand{\gd}{\mathfrak{d}}

\newcommand{\gp}{\mathfrak{p}}
\newcommand{\gq}{\mathfrak{q}}

\newcommand{\gD}{\mathfrak{D}}

\newcommand{\gR}{\mathfrak{R}}
\newcommand{\gQ}{\mathfrak{Q}}
\newcommand{\gP}{\mathfrak{P}}

\newcommand{\len}{\ell}

\newcommand{\FF}{\mathbb{F}}
\newcommand{\QQ}{\mathbb{Q}}
\newcommand{\NN}{\mathbb{N}}
\newcommand{\ZZ}{\mathbb{Z}}
\newcommand{\RR}{\mathbb{R}}
\newcommand{\CC}{\mathbb{C}}
\newcommand{\OO}{\mathcal{O}}
\newcommand{\PP}{\mathcal{P}}

\newcommand{\FFq}{\FF_q}

\newcommand{\modd}[1]{\left(\textup{mod} \ #1\right)}

\def\clap#1{\hbox to 0pt{\hss#1\hss}}

\usepackage{hyperref}

\begin{document}

\setcounter{page}{297}
\publyear{2021}
\papernumber{2100}
\volume{184}
\issue{4}

  \finalVersionForARXIV

\title{Computing the Length of Sum of Squares and Pythagoras Element in a Global Field}


\author{Mawunyo Kofi Darkey-Mensah\thanks{Address  for correspondence: Institute of Mathematics,
                          University of Silesia, ul. Bankowa 14, 40-007 Katowice, Poland. \newline \newline
          \vspace*{-6mm}{\scriptsize{Received  February 2021; \ accepted August 2021.}}},
          Beata Rothkegel
          \\
Institute of Mathematics \\
University of Silesia \\
ul. Bankowa 14, 40-007 Katowice, Poland\\
mdarkeymensah@gmail.com, beata.rothkegel@us.edu.pl
}

\maketitle

    \runninghead{M.K. Darkey-Mensah and B. Rothkegel}{Computing the Length of Sum of Squares and Pythagoras Element...}

\begin{abstract}
This paper presents algorithms for computing the length of a sum of squares and a Pythagoras element in a global field $K$ of characteristic different from $2$. In the first part of the paper, we present algorithms for computing the length in a non-dyadic and dyadic (if $K$ is a number field) completion of $K$. These two algorithms serve as subsidiary steps for computing lengths in global fields. In the second part of the paper we present a procedure for constructing an element whose length equals the Pythagoras number of a global field, termed a Pythagoras element.
\end{abstract}


\section{Introduction}
The problem of representing an element in a ring as a sum of squares is well known in mathematics, ranging from the works of Lagrange and Gauss, through the works of Waring and Hilbert, to contemporary papers (see e.g. \cite{choi1980sums} or \cite{powers1998algorithm}). For example, in papers \cite{choi1982pythagoras} and \cite{choi1995sums} the Pythagoras number is considered, while \cite{powers1998positive} gives an answer to Hilbert's  seventeenth problem, announced in 1900.

In this paper, we pursue the continuation of the recent work by P. Koprowski and A. Czogała in \cite{koprowski2018computing} on the computational aspects of the theory of quadratic forms over global fields. In \cite{koprowski2018computing} the authors focused on algorithms over number fields (i.e. finite extensions of $\QQ$). The authors in their paper developed algorithms for checking the isotropy of forms, and computing some field invariants. The aim of this article is to present algorithms for computing the length of a sum of squares and a Pythagoras element (see Definition \ref{Pythagoras_element}) in a global field of characteristic different from $2$. In this paper we make a free use of the standard results from the theory of quadratic forms over global fields. The reader is referred to \cite{lam2005introduction,scharlau1985hermitian,szymiczek1997bilinear} for a proper exposition of the theory.

Throughout this paper, if $K$ is a number field whose multiplicative group of non-zero elements is $\dot{K}$, then $\OO_K$ denotes the integral closure of $\ZZ$ in $K$, while if $K$ is a finite extension of $\FF(X)$, where $\FF$ is a finite field of characteristic not $2$, then $\OO_K$ denotes the integral closure of $\FF[X]$ in $K$. We denote by $\Omega(K)$ the set of all places of $K$. If $\gp$ is a place of $K$, then we call any valuation belonging to $\gp$ the \emph{$\gp$-adic valuation} (if $\gp$ is finite, then we call it  \emph{dyadic} and \emph{non-dyadic} if $\gp$ divides and does not divide 2, respectively). The completion of $K$ under a $\gp$-adic valuation is denoted by $K_\gp$ and called the \emph{$\gp$-adic completion} of $K$. We denote by $(\cdot,\cdot)_\gp$ the $\gp$-adic Hilbert symbol and by $h_\gp(q)$ the $\gp$-adic Hasse invariant of a quadratic form $q$ (for definitions and properties see \cite{lam2005introduction}). If $q$ is a quadratic form over $K$ (over $K_\gp$, respectively), then we write $D(q)$ ($D_\gp(q)$, respectively) for the set of all elements of $K$ ($K_\gp$, respectively) which are represented by $q$. The symbol $\langle a_1,\dots,a_n \rangle$ denotes a diagonal quadratic form over $K$ (or $K_\gp$) of dimension $n$. Next, if $\gp$ is a finite place,  then $v_\gp(a)$ denotes the $\gp$-adic valuation of an element $a$ in $K$. The square class group of the local field $K_\gp$ has the form $\dot{K}_\gp/\dot{K}^2_\gp = \{ \dot{K}^2_\gp, u_\gp \dot{K}^2_\gp, \pi_\gp \dot{K}^2_\gp, u_\gp \pi_\gp \dot{K}^2_\gp\}$, where $v_\gp (u_\gp) \equiv 0 \modd{2}$ is a $\gp$-adic unit, and $v_\gp (\pi_\gp) \equiv 1 \modd{2}$ is a $\gp$-adic uniformizer (see e.g. \cite[Theorem VI.2.2]{lam2005introduction} for further details).

Recall that the \emph{level} $s(K)$ of  $K$ is defined as the smallest positive integer $n$ (if it exists) such that $-1$ is a sum of $n$ squares of elements of $K$.

Let $a \in \dot{K}$. If $a$ is not a sum of squares of elements of $K$, then we say that $a$ has \emph{length} $\infty$ and write $\len(a) = \infty$. Otherwise, we define it's \emph{length} $\len(a)$ to be the minimal number of summands needed to express $a$ as a sum of squares of elements of $K$. If $\gp \in \Omega(K)$ is a place of $K$, then similarly  we define the \emph{length} of $a$ in the field $K_\gp$ and denote it by $\len_\gp(a)$

The \emph{Pythagoras number}  of $K$ (see e.g. \cite[XI.5.5]{lam2005introduction}), denoted $P(K)$, is the smallest positive integer $n$ such that every sum of squares in $K$ is a sum of $n$ squares. If no such integer $n$ exists, then $P(K):=\infty$. A Pythagoras element is defined as follows.

\begin{deff}\label{Pythagoras_element}
A Pythagoras element of a global field $K$, denoted $a_K$, is defined to be an element whose length is equal to the Pythagoras number of $K$. Thus $\len(a_K) = P(K)$.
\end{deff}

For example, the Pythagoras number of the rationals is $P(\QQ) = 4$ and $7 \in \QQ$ is a Pythagoras element. The Pythagoras element is not unique, e.g. $15$ is another Pythagoras element of $\QQ$.

The paper is organized as follows: in Section \ref{length_of_sos} we present algorithms (see Algorithms \ref{numberlength} and \ref{globallength}) for computing the length of a sum of squares in a number and global function field, respectively. These algorithms use subsidiary procedures (Algorithms \ref{nondyadicLength} and \ref{dyadicLength}) for deciding the lengths in a non-dyadic and dyadic completion of $K$, respectively. Next, in Section \ref{section3}, Algorithms \ref{pythagoraselement} and \ref{pythagoraselement_gff} construct a Pythagoras element in a given number field and global function field, respectively.

\section{Length of a sum of squares}\label{length_of_sos}

Let $K$ be a global field.
\begin{obs}
If $a\in \dot{K}$, then $a \in D(\langle 1,1,\dots,1 \rangle)$ for $\langle 1,1,\dots,1 \rangle$ of dimension $n$ if and only if the quadratic form $\langle a,-1,-1,\dots,-1 \rangle$ of dimension $n+1$ is isotropic.
\end{obs}
The above observation implies that if $\len(a) < \infty$, then
\begin{align}\label{l}
\begin{array}{c}
     \len(a) = \min\{n\in \NN \mid \langle a,-1,-1,\dots,-1 \rangle
     \textup{ of dimension }  \\ n+1 \textup{ is isotropic}\}
\end{array}
\end{align}
Obviously, it is true for $\len_\gp(a)$ and every $\gp\in\Omega(K)$.

\medskip
Assume $\gp \in \Omega(K)$ is a place of $K$. If $\gp$ is finite, then the $u$-invariant of $K_\gp$ is $4$ (see e.g. \cite[Theorem VI.2.12]{lam2005introduction}), so the form $\langle 1,1,1,1 \rangle$ is universal over $K_\gp$. Therefore $a \in D_\gp(\langle 1,1,1,1 \rangle)$ and $\len_\gp(a) \leq 4$. If $K$ is a number field and $\gp$ is infinite, then
\[\len_\gp(a)=
    \begin{cases}
      1 & if \ (a,-1)_\gp=1 \\
      \infty & if \ (a,-1)_\gp=-1
    \end{cases}\]
in the case when $K_\gp \cong \RR$ and  $\len_\gp(a)=1$ in the case when $K_\gp \cong \CC$.

\begin{prop}\label{prop2}
Let $K$ be a global field and $a \in \dot{K}$ with $\len(a)<\infty$, then \[\len(a) = \max_{\gp \in \Omega(K)} \len_\gp(a) \]
\end{prop}

\begin{proof}
From (\ref{l}), it follows that $\len(a)$ is the minimal natural number such that the form $\langle a,-1,-1,\dots,-1 \rangle$ of dimension $\len(a)+1$ is isotropic. By the Local-global principle \cite[Principle VI.3.1]{lam2005introduction},  $\langle a,-1,-1,\dots,-1 \rangle$ is isotropic over $K$ if and only if it is isotropic over $K_\gp$ for all $\gp \in \Omega(K)$. Therefore   $\langle a,-1,-1,\dots,-1 \rangle$ of dimension $\len(a)+1$ is isotropic over $K_\gp$ for all $\gp \in \Omega(K)$. Again from  (\ref{l}) and the Local-global principle, it follows that there is at least one $\gq\in\Omega(K)$ such that the form $\langle a,-1,-1,\dots,-1 \rangle$ of dimension $\len(a)$ is not isotropic over $K_\gq$, hence $\len_\gq(a)=\len(a)$ and $\len(a)$ is the maximal length among all $\len_\gp(a)$, $\gp\in\Omega(K)$.
\end{proof}

Now we present the first algorithm for computing the length of $a$ in a non-dyadic completion of the field $K$.

\begin{algorithm}[h] \small
\caption{Length in a non-dyadic completion} \label{nondyadicLength}
\KwIn{A nonzero element $a$ of a global field $K$ and a finite non-dyadic place $\gp$ of $K$.}
\KwOut{Length of $a$ in the completion $K_{\gp}$}
\eIf {$v_\gp(a)\equiv 1 \pmod 2$}
{
	\eIf{$-1$ is a local square in $K_\gp$}{
		\Return $2$\;
		}{
	    \Return $3$\;
	    }
}
{
\eIf{$a$ is a local square in $K_{\gp}$}{
    \Return $1$\;
    }{
    \Return $2$\;
    }
}
\end{algorithm}

\begin{myproof}[Proof of correctness of Algorithm 1:]
Assume that $v_\gp(a)$ is even. Then either $a=1$ or $a=u_\gp$ (modulo squares). If $a$ is a square, then $\len_\gp(a)=1$. Otherwise, $(-1,a)_\gp = (-1,u_\gp)_\gp = 1$, i.e. $1 \in D_\gp(\langle -1,a \rangle)$, which is equivalent to the fact that $a \in D_\gp(\langle 1,1 \rangle)$. Therefore
 $\len_\gp(a)=2$.

Now assume that $v_\gp(a)$ is odd. Then either $a=\pi_\gp$ or $u_\gp \pi_\gp$ (modulo squares) and hence
\[(-1,a)_\gp = (-1,u_\gp \pi_\gp)_\gp = (-1,\pi_\gp)_\gp =
\begin{cases}
  1 & if \ -1 \in \dot{K}_\gp^2 \\
  -1 & if \ -1 \notin \dot{K}_\gp^2
\end{cases}
\]
If $-1 \in \dot{K}_\gp^2$, then similarly as in the previous paragraph, $\len_\gp(a)=2$. If $-1 \notin \dot{K}_\gp^2$, then the level of $K_\gp$ is equal to $2$, so the form $\langle 1,1,1 \rangle$ is isotropic over $K_\gp$. Hence $a \in D_\gp(\langle 1,1,1 \rangle)$ and $\len_\gp(a)=3$.
\end{myproof}

\begin{rem}\label{isotropy_localsquare}
A procedure for testing whether an element $a$ is a square in a completion $K_\gp$ is equivalent to testing whether $x^2-a$ is irreducible in $K_\gp[x]$. Algorithms for testing the irreducibility of polynomials are already in existence and can be found for example in \cite{veres2009complexity}, \cite{guardia2011higher} or \cite{guardia2012newton}.
\end{rem}

Next, we present an algorithm for computing the length of $a$ in a dyadic completion of $K$ (if $K$ is a number field).

\begin{algorithm}[h] \small
\caption{Length in a dyadic completion} \label{alg:dyadiclength}
\KwIn{A nonzero element $a$ of a number field $K$ and a dyadic place $\gd$ of $K$.}
\KwOut{Length of $a$ in the completion $K_{\gd}$}
\If{$a$ is a local square in $K_\gd$}{
    \Return $1$\;}
Compute the Hilbert symbol $(-1,a)_\gd$\;
\If{$(-1,a)_\gd = 1$}{
    \Return $2$\;}
Compute the Hilbert symbol $(-1,-1)_\gd$\;
\If{$(-1,-1)_\gd=1$ or $-a$ is not a square in $K_\gd$}{
    \Return $3$\;}
\Return $4$\;
\label{dyadicLength}
\end{algorithm}

The proof of correctness is preceded by the following lemma.

\begin{lem}\label{lem2}
Let $a$ be a nonzero element of a number field $K$, and $\gd$ a dyadic place of $K$. The form $\langle a,-1,-1,-1 \rangle$ is isotropic over $K_\gd$ if and only if either $(-1,-1)_\gd=1$ or $-a \notin \dot{K}_\gd^2$.
\end{lem}

\begin{proof}
Assume $-a \in \dot{K}_\gd^2$. Then $\langle a,-1,-1,-1 \rangle \cong \langle -1,-1,-1,-1 \rangle$ over $K_\gd$ and $h_\gd(\langle -1,-1,-1,-1 \rangle) = (-1,-1)_\gd^6=1$. From the assumption it follows that $\langle -1,-1,-1,-1 \rangle$ is isotropic, so by the means of \cite[Proposition V.3.23]{lam2005introduction} we have $(-1,-1)_\gd = 1$.

Conversely, suppose that $(-1,-1)_\gd=1$. Then $1\in D_\gd(\langle -1, -1\rangle)$, so the  the form $\langle a,-1,-1,-1 \rangle$ is isotropic over $K_\gd$. Now assume $-a \notin \dot{K}_\gd^2$ and consider the quadratic extension $L_\gD:= K_\gd(\sqrt{-a})$ of $K_\gd$. From \cite[Example XI.2.4(7)]{lam2005introduction}, it follows that $(-1,-1)_\gD=1$ since $[L_\gD:\QQ_2] = [K_\gd(\sqrt{-a}):\QQ_2]$ is even. Moreover, $\langle a,-1,-1,-1 \rangle \cong \langle -1,-1,-1,-1 \rangle$ over $L_\gD$, and $h_\gD(\langle -1,-1,-1,-1 \rangle) = 1$. Finally, we have $h_\gD(\langle a,-1,-1,-1 \rangle)=(-1,-1)_\gD$, hence by \cite[Remark V.3.24]{lam2005introduction} $\langle a,-1,-1,-1 \rangle$ is isotropic over $K_\gd$.
\end{proof}

\begin{myproof}[Proof of correctness of Algorithm \ref{alg:dyadiclength}:]
Let $\gd$ be a dyadic place. If $a\in \dot{K}^2_\gd$, then of course $\len_\gd(a)=1$. Assume $a \notin \dot{K}^2_\gd$. We consider the $\gd$-adic Hilbert symbol $(-1,a)_\gd$. If $(-1,a)_\gd=1$, then $a\in D_\gd(\langle 1,1 \rangle)$ and $\len_\gd(a)=2$. Suppose $(-1,a)_\gd = -1$. Then $a \notin D_\gd(\langle 1,1 \rangle)$. By Lemma \ref{lem2}, if either $(-1,-1)_\gd = 1$ or $-a\notin \dot{K}_\gd^2$, then $\len_\gd(a)=3$. Otherwise, $\len_\gd(a)=4$.
\end{myproof}

\begin{rem}\label{hilbert_symbol}
An algorithm for computing the Hilbert symbol in a completion of a number field can be found in \cite[Algorithm 6.6]{voight2013identifying}.
\end{rem}

\begin{rem}\label{factorization}
In the next algorithms, Algorithms 3, 4, 5 and 6, we perform two kinds of factorization in like manner as in \cite{koprowski2018computing}. The first one is to find all dyadic primes of a given field, i.e. to factor $2\OO_K$. The second one is to find all primes dividing an element in a given field. Algorithms for factorization of ideals are well known. One may refer for example to \cite[\S 6.2.5]{cohen1993course}, \cite{cohen2012advanced} or in \cite[\S 2.2]{guardia2013new}.
\end{rem}

Now we present an algorithm for computing the length of a sum of squares in a number field.

\begin{algorithm}[h]\small
\caption{Length of a sum of squares in a number field} \label{numberlength}
\KwIn{A nonzero element $a$ of a number field $K$}
\KwOut{Length of $a$ in $K$}
\If{$K$ is formally real}{
    Let $\gR = \{\rho_1,\dots,\rho_r\}$ be the list of all real embeddings of $K$, $r \in \NN$\;
    \For {$\rho \in \gR$}{
    \If{$\rho(a)<0$}{
        \Return $\infty$}
    }
}
\If{$a$ is a square in $K$}{
    \Return $1$}
Let $\gD = \{\gd_1, \dots, \gd_m\}$ be the list of prime factors of $2$ in $\OO_K$\;
$ \len \gets 2 $\;
\For{$\gd \in \gD$}{
    Compute $\len_\gd(a)$ in $K_\gd$ using Algorithm 2\;
\If{$\len_\gd(a) = 4$}{
        \Return $4$}
   $\len\gets\max\{\len,\len_\gd(a)\}$\;
}

\If{ $\len=3$}{
\Return $3$}

Let $\gQ = \{\gq_1, \dots, \gq_n\}$ be the list of prime factors of $a$ in $\OO_K$ that do not divide $2$\;
\For{$\gq \in \gQ$}{
    Compute $\len_\gq(a)$ in $K_\gq$ using Algorithm 1\;
   \If{$\len_\gd(a) = 3$}{
        \Return $3$}
}
\Return $2$\;
\end{algorithm}

\begin{myproof}[Proof of correctness of Algorithm 3:]
Let $\rho_1,\dots,\rho_r$ be all real embeddings of $K$ for $r \geq 0$. If $\rho_i(a) < 0$ for some $i \leq r$, then $a$ is not a sum of squares in the corresponding completion, hence it cannot be a sum of squares in $K$ either.

Assume either $r=0$ or $\rho_i(a)>0$ for all $i \in \{1,\dots,r\}$. If $a \in \dot{K}^2$, then $\len(a)=1$. Therefore suppose $a\notin \dot{K}^2$.  Let $\gD$ and $\gQ $ be the set of prime factors of $2$ and  the set of prime factors of $a$ in $\OO_K$ that do not divide $2$, respectively. Moreover, put $\gP:=\gD\cup\gQ$ and fix a finite place $\gp \in \Omega(K)$. If $\gp \notin \gP$, then $v_\gp(a)=0$. It implies either $a=1$ or $a=u_\gp$ (modulo squares). If $a$ is a square, then $\len_\gp(a)=1$. Otherwise, similarly as in the proof of correctness of Algorithm \ref{nondyadicLength}, $\len_\gp(a) = 2$. Suppose $\gp \in \gP$. If $\gp$ is a dyadic place, then we use Algorithm \ref{dyadicLength}. Otherwise, we use Algorithm \ref{nondyadicLength}. Finally, by Proposition \ref{prop2}, $\len(a)=\max_{\gp \in \Omega(K)} \ \len_\gp (a)$.
\end{myproof}

Next, we present an algorithm for computing the length of a sum of squares in a global function field.
\bigskip

\begin{algorithm}[!h] \small
\caption{Length of a sum of squares in a global function field} \label{globallength}
\KwIn{A nonzero element $a$ of a global function field $K$}
\KwOut{Length of $a$ in $K$}
\If{$a$ is a square in $K$}{
    \Return $1$}
Let $\gP=\{ \gq_1, \dots, \gq_n\}$ be the list of places dividing $a$ in $K$\;
\For{$\gq \in \gP$}{
    Compute $\len_\gq(a)$ in $K_\gq$ using Algorithm 1\;
    \If{$\len_\gq(a) = 3$}{
        \Return $3$}
}
\Return $2$\;
\end{algorithm}

\begin{myproof}[Proof of correctness of Algorithm 4:]
Assume $\gP $ is the set of places dividing $a$ in $K$. The proof of correctness is similar to the proof of correctness of Algorithm~\ref{numberlength}, the second paragraph.
\end{myproof}

\section{Pythagoras element}\label{section3}

In this section, we devise algorithms that construct  Pythagoras elements in  global fields. We start with the following theorem.

\begin{thm}\label{theorem_nf}
Let $K$ be a number field, then
\begin{enumerate}[(i)]
\itemsep=0.9pt
    \item $P(K)=2$ \ iff $s(K)=1$.
    \item $P(K)=3$ \ iff $s(K) \neq 1$ and every dyadic place of $K$ has even degree.
    \item $P(K)=4$ \ iff there is a dyadic place of $K$ of odd degree.
\end{enumerate}
\end{thm}

\begin{proof}
$(i)$ Similarly as in the proofs of correctness of Algorithms \ref{nondyadicLength} and \ref{dyadicLength}, $P(K)=2$ iff  $-1 \in \dot{K}_\gp^2$ for every finite non-dyadic place $\gp$ of $K$, and $(-1,a)_\gd=1$ for every dyadic place $\gd$ of $K$ and any $a \in \dot{K}$. Hence $P(K)=2$ iff $-1 \in \dot{K}^2_\gq$ for every finite place $\gq \in \Omega(K)$ which is equivalent to $s(K)=1$.

\smallskip
$(iii)$ By Lemma \ref{lem2}, $P(K)=4$ iff $(-1,-1)_\gd=-1$ for some dyadic place $\gd$ of $K$, so $P(K)=4$ iff there is a dyadic place of $K$ of odd degree by means of \cite[Example XI.2.4(7)]{lam2005introduction}.

\smallskip
$(ii)$ Follows from $(i)$ and $(iii)$.
\end{proof}

\begin{rem}
Observe that if $K$ is a number field and $P(K)=2$, then $K$ is a nonreal field.
\end{rem}

\begin{algorithm}[h] \small
 \caption{Pythagoras element in a number field} \label{pythagoraselement}
\KwIn{A number field $K$}
\KwOut{A Pythagoras element in $K$}
\If{$s(K)=1$}{
    \Return any $a \in \dot{K} \setminus \dot{K}^2$\;
}
    Let $\gD = \{\gd_1, \dots, \gd_m\}$ be the list of prime factors of $2$ in $\OO_K$\;
    \For{$\gd \in \gD$}{
    \If{$(-1,-1)_\gd=-1$}{
    \Return $7$}
    }

     Set flag $\gets$ FALSE\;
    Set $p \gets 2$\;
    \While{flag $=$ FALSE }   {
        Let $p \gets$ smallest prime number $\geq p+1$\;
        \If{$p \equiv 3 \ (\textup{mod} \ 4)$}{
         Let $\gP = \{\gp_1, \dots, \gp_n\}$ be the list of prime factors of $p$ in $\OO_K$\;
        \For{$\gp \in \gP$}{
            \If{$v_\gp(p)\equiv 1\pmod 2$}{
                Set flag $\gets$ TRUE\;
           }
        }
         }
    }
    \Return $p$\;
\end{algorithm}

\begin{myproof}[Proof of correctness of Algorithm 5:]
Assume $P(K)=2$ and let $a \in \dot{K}\setminus \dot{K}^2$. Then $s(K)=1$ and $-1 \in \dot{K}^2$, so $-1 \in \dot{K}^2_\gq$ for every finite place $\gq \in \Omega(K)$. If $\gq$ is a finite place such that $a \in \dot{K}^2_\gq$, then $\len_\gq(a)=1$. Otherwise, similarly as in the proofs of correctness of Algorithms 1 and 2, $\len_\gq(a)=2$ so $\len(a)=2$  and $a$ is a Pythagoras element.

\medskip
Next, assume $P(K)=3$, then $s(K) \neq 1$ and $(-1,-1)_\gd=1$ for every dyadic place $\gd$ of $K$. Since $-1 \notin \dot{K}^2$, there exists a finite non-dyadic place $\gp$ of $K$ such that $-1 \notin \dot{K}^2_\gp$. If $a$ is an element of $K$ such that $v_\gp(a)$ is odd, then $a = \pi_\gp$ or $u_\gp \pi_\gp$ (modulo squares) and $\len_\gp(a)=3$. If $\gd$ is a dyadic place of $K$, then $(-1,-1)_\gd=1$ so by Lemma \ref{lem2}, $\len_\gd(a) \leq 3$. Moreover, if $a$ is a totally positive element (in the case when $K$ is formally real), then $\len(a)=3$ and $a$ is a Pythagoras element of $K$. Now let $p \equiv 3 \ (\textup{mod} \ 4)$ be any prime number factoring into prime ideals of the form $p\OO_K = \gp_1^{e_1}\cdots\gp_n^{e_n}$, and let $\gp$ be any of those factors such that $e(\gp | p) \equiv 1 \ (\textup{mod} \ 2)$. Then by \cite[Corollary VI.2.6]{lam2005introduction}, since $p \equiv 3 \modd{4}$, it implies that $-1 \notin \dot{K}_{\gp}^2$, and $p$ is the sought element.

\medskip
Finally, if $P(K)=4$, then choose a dyadic place $\gd$ of $K$ such that $(-1,-1)_\gd=-1$. We prove that $7$ is a Pythagoras element of $K$. Indeed, $-7 \in \dot{\QQ}^2_2 \subset \dot{K}^2_\gd$. By Lemma \ref{lem2}, it implies that $\langle 7,-1,-1,-1 \rangle$ is anisotropic over $K_\gd$ and $\len_\gd(7)=4$. Therefore $\len(7)=4$ and $7$ is a Pythagoras element.
\end{myproof}


\begin{thm}\label{theorem_gff}
Let $K$ be a global function field with full field of constants $\FFq$ of order $q$, then
\begin{enumerate}[(i)]
    \item $P(K)=2$ \ iff \ $q \equiv 1 \ (\textup{mod} \ 4)$
    \item $P(K)=3$ \ iff \ $q \equiv 3 \ (\textup{mod} \ 4)$
\end{enumerate}
\end{thm}

\begin{proof}
$(i)$ Similarly as in the proof of correctness of Algorithm \ref{nondyadicLength}, $P(K)=2$ if and only if $-1 \in \dot{K}^2_\gp$ for every place $\gp \in \Omega(K)$. Therefore \[P(K)=2 \iff -1 \in \dot{\FFq}^2 \subset \dot{K}^2 \iff q \equiv 1 \modd{4}\]

$(ii)$ Follows from $(i)$.
\end{proof}

\begin{algorithm}[h]\small
 \caption{Pythagoras element in a global function field} \label{pythagoraselement_gff}
\KwIn{A global function field $K$ with full field of constants $\FF_q$ of order $q$.}
\KwOut{A Pythagoras element in $K$}
\If{$q \equiv 1 \ (\textup{mod} \ 4)$}{
    \Return any $a \in \dot{\FFq} \setminus \dot{\FFq}^2$\;
}
Set flag $\gets$ FALSE\;
Set $m \gets 0$\;
\While{flag $=$ FALSE}{
    $m \gets m+1$\;
    Factor $(x^{q^m}-x)$ into monic irreducible polynomials in the form of a list: $\PP = \{p_1,\dots,p_k\}$\;
    \For{$p \in \PP$}{
    Let $\gP = \{\gp_1, \dots, \gp_n\}$ be the list of places dividing $p$ in $K$\;
    \For{$\gp \in \gP$}{
        \If{$v_\gp(p)\equiv 1 \pmod 2$}{
           Set flag $\gets$ TRUE\;
           }
        }
    }
}
\Return $p$\;
\end{algorithm}

\begin{myproof}[Proof of correctness of Algorithm 6:]
The full field $\FFq$ of constants is algebraically closed in $K$. Hence $a \notin \dot{\FFq}^2$ implies that $a \notin \dot{K}^2$. Therefore if $P(K)=2$, we have $\len(a)=2=P(K)$.

\medskip
Conversely, if $q \equiv 3 \ (\textup{mod} \ 4)$, then $-1 \notin \dot{K}^2$ by Theorem \ref{theorem_gff}, so there is a place $\gp \in \Omega(K)$ such that $-1 \notin \dot{K}_\gp^2$. If $a \in K$ and $v_\gp(a)$ is odd, then similarly as in the proof of correctness of Algorithm 5, $\len_\gp(a)=3$ which is the maximal length in $K$ and $a$ is a Pythagoras element. Further, for any positive integer $m$, it is well known that the polynomial $x^{q^m}-x$ in $\FFq[x]$ factors into a product $\prod_{d\mid m}P(d,q)$ of monic irreducible polynomials $P$ of degree $d$ (see e.g.  \cite[Chap. 7, Theorem 2]{ireland2012classical}). If $P$ is any monic irreducible polynomial factoring into a product of powers of prime ideals $\gp_1^{e_1} \cdots \gp_k^{e_k}$ in $\OO_K$ such that $e_i(\gp_i | P) \equiv 1 \ (\textup{mod} \ 2)$ for any $i\leq k$, then $(-1,P)_{\gp_i}=-1$ which implies $-1\notin \dot{K}_{\gp_i}^2$, and $P$ is the sought element.
\end{myproof}

The presented algorithms can be implemented in existing computer algebra systems. In fact, they have currently been implemented in CQF -- a free, open-source Magma package for doing computations in quadratic forms theory (see \cite{koprowski2020cqf}). CQF determines the length of an element and a Pythagoras element in a global field using the functions \texttt{LengthOfSumOfSquares} (or \texttt{LengthOfSOS} for short) and \texttt{PythagorasElement}, respectively.


\begin{thebibliography}{10}
\providecommand{\url}[1]{\texttt{#1}}
\providecommand{\urlprefix}{URL }
\expandafter\ifx\csname urlstyle\endcsname\relax
  \providecommand{\doi}[1]{doi:\discretionary{}{}{}#1}\else
  \providecommand{\doi}{doi:\discretionary{}{}{}\begingroup
  \urlstyle{rm}\Url}\fi
\providecommand{\eprint}[2][]{\url{#2}}

\bibitem{choi1980sums}
Choi MD, Lam TY, Reznick B, Rosenberg A.
\newblock Sums of Squares in Some Integral domains.
\newblock \emph{Journal of Algebra}, 1980.
\newblock \textbf{65}(1):234--256. doi:10.1016/0021-8693(80)90248-3.

\bibitem{powers1998algorithm}
Powers V, W\"{o}rmann T.
\newblock An algorithm for Sums of Squares of Real Polynomials.
\newblock \emph{Journal of Pure and Applied Algebra}, 1998.
\newblock \textbf{127}(1):99--104. doi:10.1016/S0022-4049(97)83827-3.

\bibitem{choi1982pythagoras}
Choi MD, Doi ZD, Lam TY, Reznick B.
\newblock The {Pythagoras} Number of Some Affine Algebras and Local Algebras.
\newblock \emph{Journal f\"{u}r die Reine und Angewandte Mathematik}, 1982.
\newblock \textbf{336}:45--82.    doi:eudml.org/doc/59018.

\bibitem{choi1995sums}
Choi MD, Lam TY, Reznick B.
\newblock Sums of Squares of Real Polynomials. {$K$}-theory and Algebraic
  Geometry: connections with Quadratic Forms and Division Algebras ({Santa
  Barbara, CA}, 1992).
\newblock \emph{Proceedings of Symposia in Pure Mathematics}, 1995.
\newblock \textbf{58}:103--126.  ISBN:978-0-8218-9361-6, 978-0-8218-0340-0.

\bibitem{powers1998positive}
Powers V.
\newblock Positive Polynomials and Sums of Squares: theory and Practice. {Real}
  Algebraic Geometry.
\newblock \emph{Panoramas and Syntheses}, 2017.
\newblock \textbf{51}:155--180.  ISSN:1272-3835.

\bibitem{koprowski2018computing}
Koprowski P, Czoga{\l}a A.
\newblock Computing with Quadratic Forms Over Number Fields.
\newblock \emph{Journal of Symbolic Computation}, 2018.
\newblock \textbf{89}:129--145.   doi:10.1016/j.jsc.2017.11.009.

\bibitem{lam2005introduction}
Lam TY.
\newblock Introduction to Quadratic Forms Over Fields, volume~67.
\newblock American Mathematical Soc., 2005. ISBN-10:0-8218-1095-2, 13:978-0-8218-1095-8.

\bibitem{scharlau1985hermitian}
Scharlau W.
\newblock Hermitian Forms over Global Fields.
\newblock Springer, 1985.

\bibitem{szymiczek1997bilinear}
Szymiczek K.
\newblock Bilinear Algebra: An Introduction to the Algebraic theory of
  quadratic forms, volume~7.
\newblock CRC Press, 1997.   ISBN-10:9056990764, 13:978-9056990763.

\bibitem{veres2009complexity}
Veres OE.
\newblock On the Complexity of Polynomial Factorization over p-adic Fields.
\newblock Ph.D. thesis, Concordia University, 2009.  ID:976383.

\bibitem{guardia2011higher}
Gu{\`a}rdia J, Montes J, Nart E.
\newblock Higher {Newton} Polygons in the Computation of Discriminants and
  Prime Ideal Decomposition in Number Fields.
\newblock \emph{Journal de th{\'e}orie des nombres de Bordeaux}, 2011.
\newblock \textbf{23}(3):667--696.   doi.org/10.5802/jtnb.782.

\bibitem{guardia2012newton}
Gu{\`a}rdia J, Montes J, Nart E.
\newblock Newton Polygons of Higher Order in Algebraic Number Theory.
\newblock \emph{Transactions of the American Mathematical Society}, 2012.
\newblock \textbf{364}(1):361--416.   doi:10.1090/S0002-9947-2011-05442-5.

\bibitem{voight2013identifying}
Voight J.
\newblock Identifying the Matrix Ring: Algorithms for Quaternion Algebras and
  Quadratic Forms.
\newblock Springer, 2013.  doi:10.1007/978-1-4614-7488-3\_10.

\bibitem{cohen1993course}
Cohen H.
\newblock A Course in Computational Algebraic Number Theory.
\newblock \emph{Graduate texts in Math.}, 1993.
\newblock \textbf{138}:88.   ISBN:978-3-662-02945-9.

\bibitem{cohen2012advanced}
Cohen H.
\newblock Advanced Topics in Computational Number Theory, volume 193.
\newblock Springer Science \& Business Media, 2012.
ISBN-13:978-1461264194,  10:1461264197.

\bibitem{guardia2013new}
Gu{\`a}rdia J, Montes J, Nart E.
\newblock A New Computational Approach to Ideal Theory in Number Fields.
\newblock \emph{Foundations of Computational Mathematics}, 2013.
\newblock \textbf{13}(5):729--762.   doi:10.1007/s10208-012-9137-5.

\bibitem{ireland2012classical}
Ireland K, Rosen M.
\newblock A Classical Introduction to Modern Number Theory, volume~84.
\newblock Springer-Verlag, 1990.   doi:10.1007/978-1-4757-2103-4.

\bibitem{koprowski2020cqf}
Koprowski P.
\newblock {CQF} {Magma} Package.
\newblock \emph{ACM Communications in Computer Algebra}, 2020.
\newblock \textbf{54}(2):53--56.   doi:10.1145/3427218.3427224.
\end{thebibliography}


\end{document}